\documentclass{article}
\usepackage{amsfonts,amssymb,latexsym}
\newcounter{num}[section]
\newcommand{\Num}{\refstepcounter{num}%
\textbf{\arabic{section}.\arabic{num}}}

\newcommand{\Theorem}{\textbf{Theorem~}}
\newcommand{\Proof}{{\LARGE\emph{Proof}}}
\newcommand{\Def}{\textbf{Definition~}}
\newcommand{\Conj}{\textbf{Conjecture~}}

\newcommand{\Prop}{\textbf{Proposition~}}
\newcommand{\Cor}{ \textbf{ Corollary~}}

\newcommand{\Ch}{{\mathfrak S}}
\newcommand{\Kc}{{\cal K}}
\newcommand{\Ac}{{\cal A}}
\newcommand{\Bc}{{\cal B}}
\newcommand{\Cc}{{\cal C}}

\newcommand{\al}{{\alpha}}
\newcommand{\la}{{\lambda}}

\newcommand{\Fq}{{\Bbb F}_q}
\newcommand{\Cb}{{\Bbb C}}

\newcommand{\Ad}{{\mathrm{Ad}}}

\newcommand{\rt}{{\mathrm{right}}}

\newcommand{\GL}{{\mathrm{GL}}}

\newcommand{\Ind}{{\mathrm{Ind}}}

\newcommand{\row}{{\mathrm{row}}}
\newcommand{\col}{{\mathrm{col}}}

\renewcommand{\leq}{\leqslant}

\begin{document}
\Large

\title{Supercharacter theories for algebra group extensions}
\author{A.N.Panov\footnote{The work is performed at NRU HSE with the support of the Russian Science Foundation Grant no. 16-41-01013.}}
\date{}
 \maketitle

\begin{abstract}
 We construct a few supercharacter theories for finite  semidirect products with the normal subgroup of  algebra group type. In the case of  algebra groups, these supercharacter theories coincide with the one of P.Diaconis and I.M.Isaaks.
 For  the parabolic subgroups of $\GL(n,\Fq)$,  the supercharacters and superclasses are classified.
\end{abstract}
{\small \textit{2000 Mathematics Subject Classification}. 20C33, 05E10. }

\section{Introduction}

The notion of a superharacter theory was introduced by P. Diaconis and I. M. Isaaks in the paper  \cite{DI}.  By definition, a supercharacter theory of a finite group $G$ is a pair  $(\Ch,\Kc)$, where $\Ch=\{\chi_1,\dots,\chi_m\}$  is a system of orthogonal (disjoint) characters (representations) of the group  $G$ and $\Kc=\{K_1, \ldots, K_m\}$ is a partition of  $G$ such that each character $\chi_i$ is constant on each class  $K_j$ and  $\{1\}\in \Kc$. The characters from   $\Ch$ are called \emph{supercharacters}  and the subsets of $\Kc$ \emph{superclasses}.  Remark that the number of supercharacters  is equal to the number of superclasses.

 The theory of irreducible characters is an example of a supercharacter theory. In this theory,  the supercharacters are the irreducible characters, and the superclasses are the classes of conjugate elements. In general case, super\-classes are unions of classes of conjugate element. Therefore, one can treat  supercharacter theory as an approximation of the theory of irreducible characters.

  For some groups, such as the unitriangular group  $\mathrm{UT}(n,\Fq)$,  the classification of irreducible characters is  a extremely difficult, "wild"\, problem.
  In this case,  we have a goal to construct a supercharacter theory that affords the most precise approximation of  the theory of irreducible representation.

In the above-mentioned paper   $\cite{DI}$, the authors constructed the supercharacter theory for algebra groups.  An algebra group is a group of the form  $U=1+J$, where $J$ is an associative finite dimensional nilpotent algebra over a finite field $\Fq$. Let us summarize the constructions of supercharacters and superclasses for  algebra groups.
 The superclasses of  $U$ are the subsets of the form $K(1+x) = 1+UxU$, where $x\in J$.

Turn to supercharacters.
Consider the left and right actions of  $U$  on  $J^*$ defined by the formulas
$u\la(x)=\la(xu)$ and $\la u(x)=\la(ux)$.  Given $\la\in J^*$, we consider the stabilizer  $U_{\la,\rt}$ of the right action of  $U$ on  $J^*$.  The subgroup $U_{\la,\rt}$ is an algebra group  $U_{\la,\rt}=1+J_{\la,\rt}$, where
$$J_{\la,\rt}=\{y\in J:~ \la(yx)=0~ \mbox{for~ all} ~x\in J\}.$$
 Fix the nontrivial character   $t\to \varepsilon^t$ of the additive group  $\Fq$ with values in  $\Cb^*$. Define the linear character of stabilizer  $\xi_\la (u)=\varepsilon^{\la(x)}$, where $u=1+x$.    A supercharacter of the algebra group  $U$ is  the character  $\chi_\la$ induced from the linear  character  $\xi_\la $ of  $U_{\la,\rt}$.

 In accordance with the paper \cite{DI}, the systems of characters  $\{\chi_\la\}$ and  subsets  $\{K(1+x)\}$, where $\la$ and  $x$ run through the systems of representatives of  $U\times U$ orbits in $J^*$ and $J$ respectively, afford  a supercharacter theory of the algebra group  $U$.
 In the case  $U$ is the unitriangular group, the supercharacters and superclasses are parametrized by the pairs  $(D,\phi)$, where  $D$ is a  rook placement (the other term  is a basic subset) in the set of positive roots, and $\phi$ is a map $D\to\Fq^*$ (see \cite{A1,A2,A3,Yan}).

Majority of the papers on  this topic was devoted to abelian and unipotent groups. The following issues were considered:
classification of supercharacter theories (for abelian groups), restriction and induction  in the supercharacter theory for algebra groups, applications in the number theory and in the solution of  random walk problem.  The bibliography  is contained in  \cite{VERY}.

 In the paper ~\cite{H}, there is proposed  the general construction of a super\-cha\-racter theory for semidirect products of groups. However, this construction gives rise to a coarse supercharacter theory.
Some papers are devoted to super\-cha\-racter theories of semidirect products with abelian normal subgroups (see \cite{SA,L1,L2}).

  Let  $L$ be a finite group, and  $U=1+J$ be an algebra group.  
Suppose that there are defined  the left   and  right  actions of  $L$  on $J$
such that these actions commute and  the following conditions hold:
\begin{enumerate}
\item $h(xy)=(hx)y$ ~and ~$(xy)h=x(yh)$,
\item $x(hy)=(xh)y,$
\end{enumerate}  
 for any  $x,y\in J$ and $h\in L$.  
  Then we have the homomorphism $\Ad: L\to\mathrm{Aut(U)}$ defined by the formula  $\Ad_h(1+x)=1+ hxh^{-1}$.
We refer to  the semidirect product $G=L\ltimes U$ of as an  \emph{algebra group extension}.
  Examples of such semudirect products are the parabolic subgroups of  $\mathrm{GL}(n)$ and the groups of invertible elements  of associative finite dimensional algebras. In the case   $L$ is an  abelian group and  $|L|$ does not devide  $\mathrm{char}\,\Fq$, these groups were considered in the paper  \cite{P1} and  were called finite groups of triangular type.

In this paper, we present three superclass theories for the algebra group  extensions  $G=L\ltimes U$. These  superclass theories are called  $GG$-, $GU$- and  $UU$-supercharater theories. They differs by the choice of subgroup in the stabilizer $H_\la$ of a linear form  $\la\in J^*$ (see theorems  \ref{GGthe},~\ref{GUthe},~\ref{UUthe}). In the case of parabolic subgroups of  $\GL(n,\Fq)$, the supercharacters and super\-classes of
$GG$- and $GU$-supercharacter theories are characterized in terms of rook placements  (see theorems  \ref{GGthep}, \ref{GUthep}). As for the $UU$-supercharacter theory,  classification of super\-cha\-racter and super\-clas\-ses is an open problem up today (see conjecture \ref{conjone}). Here we merely note that  for maximal parabolic subgroups the $UU$-supercha\-rac\-ter theory coincides with the theory of irreducible characters (see proposition \ref{UUprop}).

The other approach that aimed to construct  a supercharacter  theory  for such groups is presented in the papers \cite{P1,P2}. This approach leads the supercharacter theory  that is incomparable with the one proposed in this paper
(here the term incomparable means that the superclass partitions are incomparable).

\section{Formula for supercharacters}

In this section, we obtain the formula for the character defined by a linear form  $\la\in J^*$ and  a character (representation) of certain subgroup $H_0$ in the stabilizer of $\la$.

Each of the subgroups $L$ and $U$ affords the left and right actions on $J$. 
For $g=hu$, where $h\in L$ and $u\in U$, the formulas  $ g x=h(ux),~~ xg=(xh)u$
define the left and right actions of $G$ on $J$.
These actions obey conditions 1 and 2 of Introduction. 
  
Define the left and right actions of  $G$ on  $J^*$ by the formulas
\begin{equation}\label{ltrt}
g\la(x)=\la(gx)~\mbox{and}~\la g(x)=\la(xg).
\end{equation}

For an arbitrary  $\la\in J^*$, we consider its stabilizer    $H_\la=\{h\in L: ~ \Ad^*_h(\la)=\la\}$ with respect to the coadjoint representation of  $L$ in $J^*$.

For two subgroups  $G_1$ and $G_2$ in  $G$, we define the subgroup  $H_{G_1\la G_2}$ that consists of elements of  $L$
that stabilize all linear forms $\mu\in G_1\la G_2$. In the case  $G_2=\{1\}$ (respectively,  $G_1=\{1\}$), we use  notation  $H_{G_1\la}$ (respectively,  $H_{\la G_2}$).
In particular, we define the subgroups  $H_{G\la G}$, $H_{G\la U}$, ~$H_{U\la U}$ and
$H_{G\la}$,~ $H_{U\la}$. The following statement will not be used to construct a supercharacter theory, nevertheless it is interesting in itself.
\\
\Prop\Num. The subgroup  $H_{G\la U}$ coincides with  $H_{G\la}$. Respectively,  $H_{U\la U}$ coincides with $H_{U\la}$.\\
\Proof. We verify  $H_{G\la U} = H_{G\la}$. The second equality can be proved analogically.
It is obvious that   $H_{G\la U} \subseteq H_{G\la}$. Let us prove the contrary inclusion.  Let $r\in  H_{G\la}$. Then  $r(g\la)r^{-1}=g\la$ for any  $g\in G$. For any  $ 1+y\in U$,  we have
$r((1+y)g\la)r^{-1}=(1+y)g\la$. We obtain  $r(yg\la)r^{-1}=yg\la$. Therefore,  $g\la(r^{-1}xry)=g\la(xy)$ for any  $x,y\in J$. As  $Ad^*_r$  stabilizes  $g\la$, we have $g\la(xryr^{-1})=g\la(xy)$. Then  $r^{-1}g\la (1+x) r (y) = g\la (1+x)(y)$ for any  $x,y\in J$. This implies  $r(g\la u)r^{-1}=g\la u$
for any  $g\in G$  and $u\in U$. Hence $h\in H_{G\la U}$. ~$\Box$

Let  $H_0$ be a subgroup in $H_\la$. Consider the subgroup  $G_{0,\la}=H_0\ltimes U_{\la,\rt}$. Let  $T$  be a representation of the subgroup $H_0$;  we shall show that the formula   $$\Xi_{\theta,\la}(g)= T(h)\varepsilon^{\la(x)},$$
where $g=h(1+x)$, ~ $h\in H_0$ and $x\in J$, defines a representation of the subgroup
 $G_{0,\la}$.

Indeed, for  $g_1=h_1(1+x_1)$ and  $g_2=h_2(1+x_2)$, we have
 $$\Xi_{\theta,\la}(g_1g_2) = \Xi_{\theta,\la}(h_1(1+x_1)h_2(1+x_2))=
 \Xi_{\theta,\la}(h_1h_2(1+h_2^{-1}x_1h_2)(1+x_2))=$$
 $$T(h_1h_2)\varepsilon^{\la(h_2^{-1}x_1h_2)}\varepsilon^{\la(x_2)}\varepsilon^{\la(h_2^{-1}x_1h_2x_2)}.$$
 Since $x_1\in J_{\la,\rt}$ and $h_2\in H_0\subset H_\la$, we obtain
 $\la(h_2^{-1}x_1h_2)=\la(x_1)$ and $\la(h_2^{-1}x_1h_2x_2)=\la(x_1h_2x_2h_2^{-1})=0$.
 We get
 $$\Xi_{\theta,\la}(g_1g_2) = T(h_1)T(h_2)\varepsilon^{\la(x_1)}\varepsilon^{\la(x_2)}=\Xi_{\theta,\la}(g_1)\Xi_{\theta,\la}(g_2).$$

 Denote by  $\theta$ the character of representation $T$, and by  $\xi_{\theta,\la}$ the one of  $\Xi_{\theta,\la}$.
 We have
\begin{equation}\label{thla}
\xi_{\theta,\la}(g)=\theta(h)\varepsilon^{\la(x)}.
\end{equation}
Remark that  $U_{g\la,\rt} = U_{\la,\rt}$ for each $g\in G$.
If  $H_0\subseteq H_{G\la}$, then  $G_{0,g\la}=G_{0,\la}$.
Consider the character  $\chi_{\theta,\la}$ induced from the character  \begin{equation}\label{thGla}
\sum_{p\in L}\xi_{\theta,p\la}
\end{equation}
of the subgroup  $G_{0,\la}$ of $G$. \\
\Prop\Num. Suppose that  $H_0\subseteq H_{G\la}$. Then   the value of the character  $\chi_{\theta,\la}$ at
$g=h(1+x)$ is calculated by the following formula
\begin{equation}\label{chigen}
\chi_{\theta,\la}(g) =  \frac{ |L|\cdot |\la U|}{|H_0|\cdot |G\la U|} \sum_{r\in L}\dot{\theta}(rhr^{-1})\left(\sum_{ \mu\in G\la U}\varepsilon^{\mu(rxr^{-1})}\right),
\end{equation}
where $\dot{\theta}$ is a function on $L$ that  equals  to  $\theta$ on  $H_0$ and is zero outside  $H_0$.\\
\Proof.
Let $G_0=H_0U$.
Denote by $\overline{\chi}_{\theta,\la}$ the character of $G_0$ induced   from the character (\ref{thGla}) of   $G_{0,\la}$.
Then $\chi_{\theta,\la}=\Ind(\overline{\chi}_{\theta,\la},G_0,G)$.
 Let $\dot{\xi}_{\theta,\la}$ be a function on  $G_0$ that equals to $\xi_{\theta,\la}$ on $G_{0,\la}$ and is zero outside
$G_{0,\la}$.
By definition of induced character,  $$\overline{\chi}_{\theta,\la}(g_0)=\frac{1}{|G_{0,\la}|} \sum_{s_0\in G_0,~p\in L} \dot{\xi}_{\theta,p\la}(s_0g_0s_0^{-1}).$$
For $g_0=h_0u\in G_0$ and $s_0=\rho_0v\in G_0$, where $h_0, \rho_0\in H_0$ and $u, v\in U$, we have
$$s_0g_0s_0^{-1}= \rho_0v(h_0u)v^{-1}\rho_0^{-1} = \rho_0h_0\rho_0^{-1}\cdot \rho_0v^{h_0}uv^{-1}\rho_0^{-1}, $$
where $v^{h_0}=h_0^{-1}vh_0$.

The algebra  $J$ acts on the left and right on  $J^*$ by formulas
$y\la(x)=\la(xy)$ and $\la y(x)=\la(yx)$. By definition of the right stabilizer $J_{\la,\rt}$,
 we obtain  $J_{\la,\rt}^\perp = J\la$.  For each $p\in L$,
$$ \sum_{y\in J} \varepsilon^{yp\la(x)} =\left\{\begin{array}{ll} |J|&, \mbox{if} ~x\in J_{\la,\rt};\\ 0&, \mbox{if}~ x\notin J_{\la,\rt}. \end{array}\right.$$
Therefore,
$$
\sum_{a\in U} \varepsilon^{ap\la(x)} =\left\{\begin{array}{ll} |U|\varepsilon^{p\la(x)}&, \mbox{if}~ x\in J_{\la,\rt};\\ 0&, \mbox{if}~ x\notin J_{\la,\rt}. \end{array}\right.$$
Thus,
$$\overline{\chi}_{\theta,\la}(g_0)=\frac{ 1}{|G_{0,\la}|\cdot |U|} \sum_{\rho_0\in H_0,~ p\in L,~ a, v\in U}\theta(\rho_0h_0\rho_0^{-1})\varepsilon^{ap\la(\rho_0(v^{h_0}uv^{-1}-1)\rho_0^{-1})}.$$

Since $H_0\subseteq H_{G\la}$ and $h_0\in H_0$, we see
$$	
ap\la(\rho_0(v^{h_0}uv^{-1}-1)\rho_0^{-1}) = ap\la(v^{h_0}uv^{-1}-1) =
ap\la(v^{h_0}uv^{-1}-uv^{-1} + uv^{-1} -1) = $$
$$ uv^{-1}ap\la(v^{h_0}-1) + ap\la(uv^{-1} -1)=
uv^{-1}ap\la(h_0^{-1}(v -1)h_0) + ap\la(uv^{-1} -1) = $$
$$ uv^{-1}ap\la(v -1) + ap\la(uv^{-1} -1) =
ap\la((v -1)uv^{-1}+uv^{-1}-1) = ap\la(v(u-1)v^{-1}).$$
For $u=1+x$, we get
$$
\overline{\chi}_{\theta,\la}(g_0)=
\frac{ |H_0|}{|G_{0,\la}|\cdot |U|} \sum_{ p\in L,~ a, v\in U}\theta (h_0)\varepsilon^{v^{-1}ap\la v(x)}=	$$
$$
\frac{ |H_0|}{|G_{0,\la}|\cdot |U|}\cdot \theta (h_0) \sum_{ p\in L,~ a, v\in U}\varepsilon^{ap\la v(x)} = \frac{ |H_0|\cdot |G|\cdot |U|}{|G_{0,\la}|\cdot |U|\cdot |G\la U|}\cdot \theta (h_0) \sum_{ \mu\in G\la U}\varepsilon^{\mu(x)}=$$
$$
 \frac{ |L|\cdot |\la U|}{ |G\la U|}\cdot \theta (h_0) \sum_{ \mu\in G\la U}\varepsilon^{\mu(x)}.
$$
The above formula for $\overline{\chi}_{\theta,\la}$ implies the formula (\ref{chigen}) for $\chi_{\theta,\la}(g) $. ~$\Box$\\
\Prop\Num. Suppose that $H_0\subseteq H_{U\la}$. Let $\chi''_{\theta,\la}$ be the character induced from the character  $\xi_{\theta,\la}$ of the subgroup  $G_{0,\la}$. Then the value
of  $\chi''_{\theta,\la}$ at $g=h(1+x)$ is calculated by the formula
\begin{equation}\label{chiUgen}
\chi''_{\theta,\la}(g) =  \frac{|\la U|}{|H_0|\cdot |U\la U|} \sum_{r\in L}\dot{\theta}(rhr^{-1})\left(\sum_{ \mu\in U\la U}\varepsilon^{\mu(rxr^{-1})}\right).
\end{equation}
\Proof. The proof is similar to the one of the previous Proposition.

\section{$GG$-supercharacter theory}

\subsection{$GG$-supercharacter theory  for  algebra group extensions}
The subgroup $H_{G\la G}$ is a normal subgroup in  $L$.
Consider the set of pairs $\Ac=\{(\theta,\la)\}$, where $\la$ runs through the set of representatives of the $G\times G$ orbits in $J^*$, and  $\theta$ runs through the set  $\mathrm{Irr}_L(H_{G\la G})$ of all  $L$-irreducible characters of the subgroup  $H_{G\la G}$.
In the case  $H_0=H_{G\la G}$, we denote the character
$\chi_{\theta,\la}$ by   $\chi_\al$. The following formula is a corollary of (\ref{chigen}).  \\
\Cor\Num.  For $\al=(\theta,\la)\in\Ac$, the value of the character  $\chi_\al$  at $g=h(1+x)$ is calculated by the formula
\begin{equation}\label{chiGlaG}
\chi_\al(g) =  \frac{ |L|^2\cdot |\la U|}{|H_{G\la G}|\cdot |G\la G|}\cdot \dot{\theta}(h)\cdot \sum_{ \mu\in G\la G}\varepsilon^{\mu(x)}.
\end{equation}
\Proof.  In the case  $H_0=H_{G\la G}$, the  formula  (\ref{chigen}) has the form
$$
\chi_\al(g) =  \frac{ |L|\cdot |\la U|}{|H_{G\la G}|\cdot |G\la U|}\cdot\dot{\theta}(h)\cdot \sum_{r\in L}\left(\sum_{\mu\in G\la U}\varepsilon^{\mu r(x)}\right). $$
The group  $L$ transitively acts on the set of  $G\times U$ orbits in $G\la G$ by right multiplication. The number of  $G\times U$ orbits in  $G\la G$ equals to $|G\la G|/|G\la U|$. The number of elements in the stabilizer in $L$ of a  $G\la U$  orbit  is equal to
$$\frac{|L|\cdot |G\la U|}{|G\la G|}.$$
Then  $$\chi_{\theta,\la}(g) =  \frac{ |L|\cdot |\la U|}{|H_{G\la G}|\cdot |G\la U|}\cdot  \frac{|L|\cdot |G\la U|}{|G\la G|}\cdot   \dot{\theta}(h)\cdot \sum_{\mu\in G\la G}\varepsilon^{\mu (x)}. $$
Hence we obtain  (\ref{chiGlaG}).~ $\Box$\\
\Prop\Num\label{orthog}. The characters  $\{\chi_\al\}$ are pairwise orthogonal.\\
\Proof. The functions  $\{\varepsilon^{\mu(x)}\}$ form a system of irreducible characters of the  abelian group  $J$ with values in  $\Cb^*$.  The irreducible characters are pairwise disjoint. Therefore, the characters  $\chi_\al$ and $\chi_{\al'}$, which attached  to different  $G\times G$ orbits,  are orthogonal.
If they are attached to the same  $G\times G$ orbit, then the characters  $\chi_\al$ and $\chi_{\al'}$ are orthogonal, since the different  $L$-irreducible characters of the group  $H_{G\la G}$ are orthogonal. ~$\Box$

Turn to construction of  superclasses.
Let   $h\in L$. If  $h\in H_{G\la G}$ for some  $\la\in J^*$,
then $\Ad_h^*$ stabilizes  all elements of  $G\la G$ and, therefore, all elements of the subspace $<G\la G>$ spanned by  $G\la G$. Consider the sum of all subspaces  $<G\la G>$  such that  $h\in H_{G\la G}$. We obtain the maximal $G\times G$ invariant subspace on which $\Ad_h^*$ acts identically.
 Its orthogonal complement is the  $L\times L$ invariant ideal in  $J$. Denote it by $J_h$.
Define the natural  projection  $\pi_h: J\to J/J_h$.

Consider the set of pairs  $\Bc=\{(h,\omega)\}$, where  $h$ runs through the set of representatives of classes of conjugate elements in  $L$, and  $\omega$  runs through the set of  $G\times G$ orbits  in $J/J_h$. For $\beta=(h,\omega)\in \Bc$, we define  the subset
\begin{equation}\label{Kclass}
K_\beta= Cl_L(h)(1+ \pi_h^{-1}(\omega)),
\end{equation}
 where $Cl_L(h)$ is the class of conjugate elements of $h$ in $L$. The system of subsets $\{K_\beta\}$ form a partition of $G$.
\\
\Theorem\Num\label{GGthe}. The systems of characters  $\{\chi_\alpha:~ \alpha \in {\mathcal{A}}\}$ and  subsets  $\{K_\beta:~ \beta \in {\mathcal{B}}\}$ give rise to a supecharacter theory of the group  $G=L\ltimes U$.\\
This supercharacter theory will be referred to as the $GG$-supercharacter theory of the group $G=L\ltimes U$.\\
\Proof. It follows from proposition \ref{orthog} that the characters  $\{\chi_\al\}$  are pairwise orthogonal. The subset  $\{1\}$ coincides with  $K_\beta$ for $\beta =(h,\omega)$, where $h=1$,~ $\omega=\{0\}$.

Let us show that the characters  $\chi_\al$ are constant on subsets $K_\beta$.
Let  $\al=(\theta,\la)$ and $\beta=(h,\omega)$.
If $h\notin H_{G\la G}$, then $\dot{\theta}(Cl_L(h))=0$ and $\chi_\al(K_\beta)=0$.
 If $h\in H_{G\la G}$, then $G\la G(J_h)=0$. It implies from the formula (\ref{chiGlaG}) that $$\chi_\al(K_\beta)= \frac{ |L|^2\cdot |\la U|}{|H_0|\cdot |G\la G|}\cdot \theta(Cl_L(h))\cdot \sum_{ \mu\in G\la G}\varepsilon^{\mu(\omega)} = const.$$

It remains to show that   $|\Ac|=|\Bc|$. Consider the subset
$$\Cc =\{(h,\la):~ h\in H_{G\la G}\},$$ where  $h$ runs through the set of representatives of classes of conjugate elements in $L$, and $\la$ runs through the set of representatives of $G\times G$ orbits in  $J^*$.

On the one hand, the number of element in  $\Cc$ equals to
$$|\Cc|=\sum_{G\la G\in G\backslash J^*/G} |Cl_L(H_{G\la G})| =
\sum_{G\la G\in G\backslash J^*/G} |\mathrm{Irr}_L(H_{G\la G})|=|\Ac|.$$
On the other hand,
$$|\Cc| =\sum \mathrm{card}\{G\la G:~ G\la G (J_h)=0\},$$
where $h$ runs through the set of representatives of classes of conjugate elements in  $L$.  The number of $G\times G$ orbits in  $(J/J_h)^*$ is equal to the number of  $G\times G$ orbits in  $J/J_h$, since, for any finite group, the number of orbits in the linear space over a finite field equals to the number of orbits in the dual space. Therefore,
$$|\Cc| =\sum \mathrm{card}\{\mbox{the~set~of}~G\times G ~ \mbox{orbits~in}~ J/J_h\}=|\Bc|.$$
We conclude $|\Ac|=|\Bc|$.~ $\Box$

\subsection{$GG$-supercharacter theory for parabolic subgroups in $\GL(n,\Fq)$}

In this subsection, we specify the $GG$-supercharacter theory for parabolic  subgroups in the group $\GL(n)$ over  $\Fq$.
 Let  $G$ be a parabolic subgroup constructed by decomposition of integer segment $[1,n]$ into union of consecutive segments  $[1,n]=I_1\cup\ldots\cup I_\ell$. The subgroup $G$ is a semidirect product  $G=LU$; here $$ L=\GL(n_1)\times\cdots \GL(n_\ell),$$  where $n_k=|I_k|$ for all  $1\leq k\leq \ell$.
The subgroup $U=\mathrm{Rad}(G)$ is an algebra group  $U=1+J$.
We refer to any pair $\gamma=(i,j)$,~ $1\leq i<j\leq n$ as a  \emph{ $J$-root}
if the corresponding matrix unit  $E_{ij}$ belongs to  $J$.
 Here $i=\row(\gamma)$  is a  \emph{row number} of $\gamma$, and  $j=\col(\gamma)$ is a \emph{column number} of $\gamma$.
 Denote the set of all  $J$-roots by
$\Delta_J$.  The system of matrix units $\{E_\gamma:~\gamma\in \Delta_J\}$ is a basis in  $J$. Respectively,  $\{E_\gamma^*\}$ is a dual basis in  $J^*$.\\
\Def\Num. A subset  $D$ in  $\Delta_J$ is called a rook placement if there is at most one root from  $D$ in any row and any column.

 Given a  rook placement  $D$, we consider the elements
 $$ x_D=\sum_{\gamma\in D} E_\gamma\in J,~~~~ u_D=1+x_D\in U,~~~~ \la_D=
 \sum_{\gamma\in D} E^*_\gamma \in J^*. $$
The set  $\Delta_J$ decomposes into a union of  $ I_k\times I_m$, ~$1\leq k<m\leq \ell$. The rook placement  $D$  also decomposes into a union of  $D_{km}=D\cap (I_k\times I_m)$. Denote  $d_{km}=|D_{km}|$.\\
\Prop\Num~\cite{Thiem}. 1) Each  $G\times G$ orbit in $J$ (respectively,  $J^*$) contains the element  $x_D$ (respectively, $\la_D$) for some rook placement $D$  in  $\Delta_J$. \\
2) Two elements  $x_D$ and $x_{D'}$ (respectively, $\la_D$ and $\la_{D'}$) belong to a common $G\times G$ orbit if and only if  $d_{km}=d'_{km}$ for any  $1\leq k < m\leq \ell$.
We say that such rook placements are  equivalent.

Let us describe  the subgroup  $H(D)=H_{G\la_D G}$.
For each  $\gamma=(i,j)\in \Delta_J$, where $i\in I_k$ and $j\in I_m$,  we consider the subgroup
$H(\gamma)$ that consists of all $h=(h_1,\ldots,h_\ell)$ such that $(h_k,\ldots,h_m)$ is a scalar submatrix and for all $t\notin [k,m]$ the component  $h_t$ is an arbitrary submatrix in $\GL(n_t)$.
 Then  $$H(D)=\bigcap_{\gamma\in D} H(\gamma).$$
{\bf Example 1. } $G$ is the parabolic subgroup of type  $(2,2,2)$ in $\GL(6)$. Let  $D=\{(1,6)\}$. The subgroup  $H(D)$ consists of scalar matrices. \\
{\bf Example 2. } $G$ is the parabolic subgroup of type  $(1,1,2,1,1)$ in
$\GL(6)$. Let $D=\{(1,2), (5,6)\}$. The subgroup  $H(D)$ is the group of block-diagonal matrices of the form $\mathrm{diag}(a,a,B,c,c)$, where $a,c\in \Fq^*$ and $B\in \GL(2)$.

The set of supercharacter parameters  $\Ac$ can be characterized as
$\Ac=\{(\theta, D)\}$, where $D$ runs through the set of representatives of rook placement classes in $\Delta_J$, and  $\theta$ runs through the set of all irreducible characters (representations) of the subgroup $H(D)$. The character $\chi_\al$ is defined as in the previous subsection.

Let us characterize superclasses.
For any  $h=\mathrm{diag}(h_1,\ldots, h_\ell)$ there exists a composition of  $[1,\ell]$ into consecutive segments  $[1,\ell]=\cup [k_i,m_i]$ such that for any segment
$[k,m]$ either  $k=m$, or $k<m$ the submatrix  $(h_k,\ldots,h_m)$ is scalar. Suppose that this decomposition is maximal among all decompositions that obey the above-mentioned conditions

 The decomposition $[1,\ell]=\cup[k_i,m_i]$ induces the decomposition of the segment $[1,n]$ with components of type  $I_{k_i}\cup\ldots\cup I_{m_i}$. The corresponding parabolic subgroup has the radical $1+\widetilde{J}$.
Observe that $\widetilde{J}$ coincides with  $J_h$.
 Denote by  $\Delta_h$ the set of roots  $\gamma\in \Delta_J$ such that  $E_\gamma\notin J_h$.

 The set of superclasses parameters  $\Bc$ can be characterized as  $\{(h,D)\}$, where  $h\in L$ runs through the set of representatives of classes of conjugate elements, and $D$ runs through the set  of representatives of rook placement classes in $\Delta_h$.
For each pair $\beta=(h,D)$, the superclass $K_\beta$ is defined by the formula (\ref{Kclass}).
\\
{\bf Example 3}.  $G$ is the parabolic subgroup of type  $(1,1,2,1,1)$ in $\GL(6)$. Let  $h=\mathrm{diag}(a,a,B,a,a)$ and
 $D=\{(1,2), (5,6))\}$. Then
 {\large
   $$  K_\beta=\left\{\left(\begin{array}{cccccc}
 a&\times&*&*&*&*\\
  0&a&*&*&*&*\\
  0&0&b_{11}'&b_{12}'&*&*\\
  0&0&b_{21}'&b_{22}'&*&*\\
   0&0&0&0&a&\times\\
    0&0&0&0&0&a\\
 \end{array}\right)\right\}\quad \mbox{for}~~ B\ne aE,$$
 $$ K_\beta=\left\{\left(\begin{array}{cccccc}
 	a&\times&*&*&*&*\\
 	0&a&0&0&0&*\\
 	0&0&a&0&0&*\\
 	0&0&0&a&0&*\\
 	0&0&0&0&a&\times\\
 	0&0&0&0&0&a\\
 \end{array}\right)\right\}\quad \mbox{for}~~ B = aE,$$}
where the submatrix  $B'=(b'_{ij})$ runs through the conjugacy class of $B$, the sigh $\times$ denotes the elements of $\Fq^*$, and the sign  $*$ denotes arbitrary elements of  $\Fq$.
  Theorem  \ref{GGthe} has the form.\\
\Theorem\Num\label{GGthep}. The systems of characters  $\{\chi_\al:~ \alpha =(\theta,D)\in \Ac\}$ and subsets  $\{K_\beta:~ \beta=(h,D) \in \Bc\}$ give rise to a supercharacter theory of the parabolic subgroup  $G=L\ltimes U$.\\

\section{$GU$-supercharacter theory}

\subsection{$GU$-supercharacter theory for algebra groups}

For an arbitrary $\la\in J^*$, we consider the stabilizer  $S_{G\la U}$ of  $G\la U$.
The subgroup  $H_{G\la U}$ (which coincides with $H_{G\la }$) is a normal subgroup in $S_{G\la U}$.
Consider the subset   $\Ac'_0=\{(\theta,\la)\}$, where  $\la$ runs through the set of representatives of  $G\times U$ orbits in $J^*$, and  $\theta$ runs through the set of  $S_{G\la U}$-irreducible characters of the subgroup  $H_{G\la U}$. The subgroup  $L$  acts on   $\Ac'_0$ conjugation. We denote by  $\Ac'$ the corresponding quotient set.

For $\al=(\theta,\la)\in \Ac'$, we define the character   $\chi'_\al=\chi_{\theta,\la}$ induced from the character
 (\ref{thGla}) of the subgroup $G_{0,\la}$, where $H_0=H_{G\la U}$. The formula (\ref{chigen}) is valid for the character  $\chi'_\al$.\\
 \Prop\Num\label{orthogp}. The characters  $\{\chi'_\al\}$ are pairwise orthogonal.\\
 \Proof. Consider the system of functions  $f_{\theta,\mu}(g) = \dot{\theta}(h)\varepsilon^{\mu(x)}$, where $g=h(1+x)$, and $\dot{\theta}(h)=\theta(h) $ for $h\in H_{G\la U}$ and $\dot{\theta}(h)=0 $ for
$h\notin H_{G\la U}$.

The functions   $\{f_{\theta,\mu}\}$ are pairwise orthogonal.
Indeed, $$(f_{\theta_1,\mu_1}, f_{\theta_2,\mu_2})=(\dot{\theta_1},\dot{\theta_2})\cdot\sum_{ x\in J} \varepsilon^{\mu_1-\mu_2(x)}.$$
The last sum equals to zero if $\mu_1\ne\mu_2$. If $\mu_1=\mu_2$, then $H_{G\mu_1U}=H_{G\mu_2U}$, and if  $\theta_1\ne\theta_2$, then  $\theta_1$ are $\theta_2$ orthogonal. The characters  $\{\chi'_\al\}$ are pairwise disjoint because they are the sums of pairwise orthogonal functions (see (\ref{chigen})).~ $\Box$

Analogically to the precious section, we construct a partition into superclasses.
Let $h\in L$. If $h\in H_{G\la U}$ for some  $\la\in J^*$,
then  $\Ad_h^*$  stabilizes all elements  of  $G\la U$, and, therefore, it stabilizes all elements in the subspace
 $<G\la U>$ spanned by  $G\la U$.  Form a sum   of all subspaces   $<G\la U>$ such that  $h\in H_{G\la U}$. We obtain the maximal  $G\times U$ invariant subspace on which   $\Ad_h^*$ acts identically.
The orthogonal complement of it is a two-sided ideal in  $J$ invariant with respect to right-side multiplication by  $L$. Denote it by  $J'_h$.
Let  $\pi'_h$ be  the natural projection  $J\to J/J'_h$.

Define the set of pairs  $\Bc'_0=\{(h,\omega)\}$, where $h\in L$, and  $\omega$  runs through the set of  $U\times G$ orbits in  $J/J'_h$. The subgroup  $L$ acts on  $\Bc'_0$ by conjugation. Denote by $\Bc'$ the corresponding quotient set.

For $\beta=(h,\omega)\in \Bc'$, we consider the subset
\begin{equation}\label{Kbeta}
K'_\beta = \bigcup_{p \in L} p (h(1+ (\pi'_h)^{-1}(\omega))p^{-1}.
\end{equation}
The system of subsets $\{K'_\beta\}$ form a partition of $G$.\\
\Theorem\Num\label{GUthe}. The systems of characters  $\{\chi'_\alpha:~ \alpha \in \Ac'\}$ and subsets  $\{K'_\beta:~ \beta \in \Bc\}$ form a supercharacter theory of the group  $G=L\ltimes U$.\\
This supercharacter theory will be referred to as the $GU$-supercharacter theory of the group $G=L\ltimes U$.\\
\Proof. The proposition  \ref{orthogp}  implies that  the characters  $\{\chi'_\al\}$  are pairwise orthogonal. The subset  $\{1\}$ coincides with $K'_\beta$ for $\beta =(h,\omega)$, where $h=1$,~ $\omega=\{0\}$.

The character $\chi'_\al$ is constant on subsets  $K'_\beta$, since the function  $$\theta(h)\sum_{\mu\in G\la U}\varepsilon^{\mu(x)}$$ is constant on  $h(1+ (\pi'_h)^{-1}(\omega))$.

It remains to show that   $|\Ac'|=|\Bc'|$. Consider the subset
$\Cc_0' =\{(h,\la)\},$ where  $\la$ runs through the set of representatives of  $G\times U$ orbits in  $J^*$ such that $h\in H_{G\la U}$. The subgroup  $L$  acts on   $\Cc'_0$ by conjugation. Denote by $\Cc'$ the corresponding  quotient set.

On the one hand,  the number of elements of  $\Cc'$ equals to
$$ |\Cc'|=\sum \{ \mbox{the ~number~of}~ S_{G\la U}-\mbox{orbits~ in }~H_{G\la U}\}=$$
$$ \sum \{ \mbox{the ~number~of}~ S_{G\la U}-\mbox{irreducible~representations~of}~H_{G\la U}\} = |\Ac'|,$$
where the sum is taken over all $G\la U\in G\backslash J^*/U$.

Let $\mathrm{Stab}_h = \{r\in L:~ rhr^{-1} = h\}$.
 The subgroup $\mathrm{Stab}_h\ltimes (G\times U)$ acts on  $J^*$ by the formula $\la\to r(g\la u)r^{-1}$. By (\ref{ltrt}) this action is comparable with the action  $x\to r(u\la g^{-1})r^{-1}$ of the same group on $J$.
We obtain
$$|\Cc'| =\sum \{\mbox{the~number~of}~ \mathrm{Stab}_h\ltimes (G\times U)-\mbox{orbits~in }~(J/J'_h)^*\}=$$
$$
\sum \{\mbox{the~number~of}~ \mathrm{Stab}_h\ltimes (G\times U)-\mbox{orbits~in}~J/J'_h\}$$
where the sum is taken over all representatives  $h$ of classes of conjugate elements in $L$. Then $|\Cc'|=|\Bc'|$.
We conclude  $|\Ac'|=|\Bc'|$.~ $\Box$

\subsection{$GU$-supercharacter theory of the parabolic subgroups in  $\GL(n,\Fq)$}

As above,  a rook placement $D$ afford the elements  $x_D$,~ $u_D$, ~ $\la_D$. The subgroup $H'(D)=H_{G\la U}$ is a stabilizer of all elements in  $G\la U$, and $S'(D)=S_{G\la_D U}$ is the stabilizer of subset  $G\la_D U$.

The set of supercharacter parameters  $\Ac'$  is identified with the set of pairs
$\{(\theta,D)\}$, where $D$  runs through the set of representatives of rook placement classes, and $\theta$ runs through the set of all  $S(D)$-irreducible characters of the subgroup  $H'(D)$.

Let us describe the subgroup  $H'(D)$.
For each $\gamma=(i,j)\in \Delta_J$, where $i\in I_k$ and $j\in I_m$,  we consider the subgroup
$H'(\gamma)$ of all elements  $h=(h_1,\ldots, h_\ell)\in L$ such that\\
1) $(h_{k+1},\ldots,h_m)$ is a scalar submatrix   $aE$;\\
2) there is a unique nonzero element in the  $i$th row of  $h_k$ which equals to $a$ and is placed  on the diagonal; \\
3) for each  $t\notin [k,m]$,   the component  $h_t$ is an arbitrary element in $\GL(n_t)$.
Then
$$H'(D)=\bigcap_{\gamma\in D} H'(\gamma).$$
{\bf Example 1. } $G$ is the parabolic subgroup of type  $(2,2,2)$ in $\GL(6)$. Let  $D=\{(1,6)\}$. Then
{\large
 $$H'(D)=\left\{\left(\begin{array}{cccccc}
a&0&0&0&0&0\\
\,*&\times &0&0&0&0\\
&&a&0&0&0\\
&&0&a&0&0\\
&&&&a&0\\
&&&&0&a
\end{array}\right)\right\},~ S'(D)=\left\{\left(\begin{array}{cccccc}
\times&0&0&0&0&0\\
\,*&\times &0&0&0&0\\
&&*&*&0&0\\
&&*&*&0&0\\
&&&&\times&0\\
&&&&*&\times
\end{array}\right)\right\}.$$}

 Let us charaterize superclasses.
  For each $h=(h_1,\ldots,h_\ell)$ there exists a decomposition of  $[1,\ell]$ into consecutive segments  $[1,\ell]=\cup[k_i,m_i]$ and  $r\in L$ such that for each segment  $[k,m]$  either   $k=m$, or $k<m$ and
\begin{equation}\label{rowGU}
(\Ad_r(h_k),\ldots,\Ad_r(h_m)) = \mathrm{diag}\left(\left(\begin{array}{cc}
aE_{s_k}&0\\
Y&X\end{array}\right), aE_{k+1}, \ldots, aE_m\right),
\end{equation}
 where $a\in \Fq^*$,~ $s_k\leq n_k$,~ $E_{s_k},\ldots, E_m$ are the unit matrices of corresponding sizes, and the submatrix  $(Y|X-aE)$ has maximal rank (equal to  $n_k-s_k$). Suppose that this decomposition is maximal among all decompositions that obey the above-mentioned conditions
 Two block-diagonal matrices with blocks of the type  (\ref{rowGU}) are conjugated in the subgroup  $L$ if and only if they are conjugated within each block by submatrices
\begin{equation}
\mathrm{diag}\left(\left(\begin{array}{cc}
\GL(s_k,\Fq)&0\\
\,*&*\end{array}\right), \GL(n_{k+1}), \ldots, \GL(n_m)\right).
\end{equation}
We shall give a description of the ideal $J'_h$ for an element $h$ with blocks of type (\ref{rowGU}).
The decomposition $[1,\ell]=\cup[k_i,m_i]$ induces the decomposition of the segment  $[1,n]$ with the components $I_{k_i}\cup\ldots\cup I_{m_i}$. The corresponding parabolic subgroup has the radical  $1+\widetilde{J}$.

 The ideal $J_h'$ coincides with  $\mathrm{Im}(\Ad_h-1)+\widetilde{J}$. It follows that $J_h'$ is spanned root subspaces. Denote by $\Delta'_h$ the set of all roots  $\gamma\in \Delta_J$ such that  $E_\gamma\notin J_h'$.

    The  $\Bc'$ can be realized as the set or pairs $\Bc'=\{(h,D)\}$, where  $h\in L$ runs through the set of
    representatives of classes of conjugate elements with blocks  (\ref{rowGU}), and  $D$ runs through the set of representatives of rook placement classes   in  $\Delta'_h$.   The superclass  $K'_\beta$ is defined as in  (\ref{Kbeta}).\\
  {\bf Example 2}. Let  $n=6$ and $G$  is a parabolic subgroup of type  $(2,1,2,1)$. Let $D=(1,3)$ and
  \begin{equation}
 h= \mathrm{diag}\left(\left(\begin{array}{cc}
  a&0\\
 y &x\end{array}\right), a, \left(\begin{array}{cc}
  a&0\\
  0&a\end{array}\right) , b \right),
  \end{equation}
  where $a,b,x\ne 0$, ~ $a\ne b$,~ $(y,x-a)\ne (0,0)$. Then, for $\beta=(h,D)$,~ $K_\beta'$ is defined by  (\ref{Kbeta}), where
  {\large
  $$h(1+\pi'_h(\omega))=\left\{\left(\begin{array}{cccccc}
  a&0&\times&*&*&*\\
  y&x &*&*&*&*\\
  &&a&0&0&*\\
  &&&a&0&*\\
  &&&0&a&*\\
  &&&&&b
  \end{array}\right)\right\}.$$ }
  \Theorem\Num\label{GUthep}. The systems of characters  $\{\chi_\al:~ \alpha =(\theta,D)\in \Ac'\}$ and subsets  $\{K_\beta:~ \beta=(h,D) \in \Bc'\}$ give rise to a supercharacter theory for the parabolic subgroup  $G=L\ltimes U$.

\section{$UU$-supercharacter theory}

\subsection{$UU$-supercharacter theory for algebra groups}

As in the previous section, we consider the stabilizer $S_{U\la U}$ of   $U\la U$.
The subgroup  $H_{U\la U}$ (with coincides with  $H_{U\la}$) is a normal subgroup in $S_{U\la U}$.
Define the subset   $\Ac''_0=\{(\theta,\la)\}$, where $\la$ runs the set of representatives of $U\times U$ orbits in $J^*$, and  $\theta$ runs through  the set of $S_{U\la U}$-irreducible characters of the subgroup  $H_{U\la U}$. The subgroup  $L$  acts on   $\Ac''_0$ by conjugation. Denote by  $\Ac''$ the corresponding quotient set.

For $\al=(\theta,\la)$, we consider the character   $\chi''_\al=\chi_{\theta,\la}$
induced from the character  (\ref{thla}) of the subgroup  $G_{0,\la}$, where  $H_0=H_{U\la U}$.

As above, each  $h\in H$ afforded  the least two-sided ideal  $J''_h$ such that  $\Ad^*_h$  acts identically on its orthogonal complement. Then  $\Ad_h$ acts identically on $J/J''_h$.
Let  $\pi''_h$ be natural projection $J\to J/J''_h$.

Define the set of pairs  $\Bc''_0=\{(h,\omega)\}$, where $h\in L$, and $\omega$  runs through the set of $U\times U$ orbits in  $J/J''_h$. The group  $L$  acts on   $\Bc''_0$ by cojugation. We denote by $\Bc''$ its quotient set.

For $\beta=(h,\omega)\in \Bc$, the class $K''_\beta$ is defined analogically to  (\ref{Kbeta}).
\\
\Theorem\Num\label{UUthe}. The systems of characters  $\{\chi''_\alpha:~ \alpha \in \Ac''\}$ and subsets  $\{K''_\beta:~ \beta \in \Bc''\}$ give rise to a supercharacter theory for the algebra group extension  $G=L\ltimes U$.\\
This supercharacter theory will be referred to as the $UU$-supercharacter theory of the group $G=L\ltimes U$.\\
\Proof. Similarly to the proof of theorem  \ref{GUthe}.

\subsection{$UU$-supercharacter theory the parabolic subgroups in  $\GL(n,\Fq)$}

The theorem \ref{UUthe} is true for parabolic subgroups, nevertheless is remain open problem on classification of supercharacters and superclasses. This problem reduces to classification of elements in   $J$ (respectively, in  $J^*$) with respect to  equivalence relation: $x\sim x'$ if there exist the elements  $a,b\in U$ and $r\in L$ such that $x'=raxbr^{-1}$.
In the paper \cite{P2}, the following conjecture was presented.\\
\Conj\Num\label{conjone}. 1) For any  $x\in J$ there exists an equivalent element $x_D$ for some rook placement in  $\Delta_J$. \\
2) Two elements  $x_D$ and $x_{D'}$ are equivalent if and only if  $D$ and $D'$ are conjugated with respect to the Weyl group of  $L$.

In the same paper  \cite{P2}, this conjecture is proved for parabolic subgroups  with blocks of size less or equal to 2.  We  state below an evident fact.\\
\Prop\Num\label{UUprop}. If the number of blocks equals  to two, then $UU$-supercharacter theory coincides with the theory of irreducible representations.

MATHEMATICAL DEPARTMENT, SAMARA UNIVERSITY, SAMARA,  RUSSIA\\
\textit{E-mail address} : \textbf{apanov@list.ru}

\end{document}